\documentclass[11pt]{amsart}
\usepackage{amsmath}

\title{Boundedness  results for periodic points on algebraic varieties}
\author{Najmuddin Fakhruddin}
\email{naf@math.tifr.res.in}
\address{School of Mathematics, Tata Institute of Fundamental Research, Homi Bhabha Road,
Mumbai 400005, India}

\newtheorem{thm}{Theorem}
\newtheorem{lem}{Lemma}
\newtheorem{prop}{Proposition}

\newtheorem{conj}{Conjecture}
\newtheorem{ques}{Question}

\theoremstyle{remark}
\newtheorem*{rem}{Remark}
\newtheorem*{rems}{Remarks}

\begin{document}
\maketitle

Let $X$ be an algebraic variety over a field $K$ and let $f:X \to X$ be
a morphism. A point $P$ in $X(K)$ is \emph{f-periodic}   if $f^n(P) = P$
for some $n>0$, and the smallest such $n$ is called the \emph{period} of
$P$. We shall prove that if $X$ and $f$ satisfy certain hypotheses, then
the set of possible periods is finite.

Our results may be viewed as an analogue of the finiteness of the torsion
of abelian varieties over finitely generated fields. It is then natural
to ask for an analogue of the full Mordell-Weil theorem. We believe that
the following conjecture is the the appropriate generalization.

\begin{conj}
Let $X$ be a proper algebraic variety over a finitely generated field $K$ of
characteristic zero and $f:X \to X$ a morphism. Suppose there exists
a subset $S$ of $X(K)$ whcih is Zariski dense in $X$ and such that
$f$ induces a bijection of $S$ onto itself. Then $f$ is an automorphism.
\end{conj}

This can be easily checked for $X = \mathbf{P}^n$ or $X$ an abelian
variety using heights and the Mordell-Weil theorem respectively.

\section{Fintely generated fields}

\begin{thm} Let $X$ be a proper variety over a field $K$ which is
finitely generated over the prime field and let $f:X \to X$ be
a morphism. \\
(i)If char$(K) = 0$ then the set of periods  of all
f-periodic points in $X(K)$ is finite.  \\
(ii)If char$(K) = p \neq 0$
then the prime to $p$ parts of the set of periods is finite i.e. there
exists $n>0$ such that all the $f^n$-periodic points in $X(K)$ have 
periods which are powers of $p$. \label{thm:one}
\end{thm}

Many special cases of this result have been known for a long time,
the first such being the theorem of Northcott \cite[Theorem 3]{northcott},
proving the finiteness of the number of periodic points in certain cases.
We refer the reader to \cite{morton-silverman}
for a more detailed list of references.

\begin{rem}We do not know whether  the periods can really
be unbounded  if $char(K)>0$.
\end{rem}

The theorem is obvious if $K$ is a finite field and we will reduce the
general case to this one by a specialization argument. A little thought
shows that the following proposition suffices to prove the theorem.

\begin{prop} Let $R$ be a discrete valuation ring with quotient field
$K$ and residue field $k$. Let $\mathcal{X}$ be a proper scheme 
of finite type over
$Spec(R)$ and $f:{\mathcal{X}} \to {\mathcal{X}}$ an R-morphism.
Assume that the conclusions of the theorem hold for $f$ restricted to
the special fibre, and that for each  $n>0$ there are only finitely many
roots of unity contained in all extensions of $k$ of degree $\leq n$.
Then the same  holds for $f$ restricted to the generic fibre,
except possibly in the case $char(K) = 0$ and $char(k)=p > 0$, when the
result holds  modulo powers of $p$. \label{prop:one}
\end{prop}

\begin{proof} Let $p = 1$ if $char(k) = 0$ .
The hypotheses imply that by replacing $f$  with a suitable
power we may assume that all the $f$-periodic points in $\mathcal{X}(k)$
have period a power of $p$. Let $P$ be a $f$-periodic point in $\mathcal{X}(K)$.
By replacing $f$ by $f^{p^n}$, for some $n$ which may depend on $P$, we may
assume that the specialization of $P$ in $\mathcal{X}(k)$ (which exists
since $\mathcal{X}$ is proper)
is a fixed point of $f$ restricted to the special fibre. Let $\mathcal{Z}$
be the Zariski closure (with reduced scheme structure) of the $f$-orbit
of $P$. $\mathcal{Z}$ is finite over $Spec(R)$ with a
unique closed point, hence is equal to
$Spec(A)$ where $A$ is a finite, local $R$-algebra (with rank equal
to  the period of $P$) which, since
$A$ is reduced,
is torsion free as an $R$-module.

The key observation of the proof is that $f$ restricted to $\mathcal{Z}$ induces an
\emph{automorphism} of finite order of $A$ (which we also denote by $f$):
Since $f$ preserves the orbit of $P$ and $\mathcal{Z}$ is reduced, it follows
that $f$ induces a map from $\mathcal{Z}$ to itself, hence an endomorphism
of $A$. $f^n$ is the identity on the orbit of $P$ for some $n>0$, hence
$f^n$ is the identity on $A \otimes_R K$. Since $A$ is torsion free,
 it follows that $f^n$ is the identity on $A$ as well.

 Let $m$ be
the maximal ideal of $A$.
Since $\mathcal{Z} $ is a closed subscheme of $\mathcal{X}$, it follows
that the dimension of $m/m^2$ is bounded independently of $P$. By the hypothesis
on roots of unity, we may replace $f$ by some power, independently of $P$, so
that the endomorphism of  $m/m^2$ induced by $f$  is the identity.
Thus, $f$ is a unipotent map with respect to the (exhaustive) filtration of
$A$ induced by powers of $m$. This implies that  the order of $f$, hence
the period of $P$, is a power of $p$.
\end{proof}

\begin{rems}
(1) For any explicitly given example, the  proof furnishes an effective method for
computing  a bound
for the periods.
(2) In the non-proper case one can prove the following result by the
methods of this paper: Let $S$ be a flat, separated, integral scheme of finite type
over $\mathbf{Z}$,
 let $\mathcal{X}$ be a separated scheme of finite type over $ S$ and let
$f: \mathcal{X} \to \mathcal{X}$ be an $S$-morphism.
If one defines the notion of $f$-periodic points and periods for elements of
$\mathcal{X}(S)$ in the obvious way, then the set of periods is again
bounded.
(3) One may also ask whether Theorem \ref{thm:one} itself holds without the
assumption of properness, for example when $X$ is arbitrary but $f$ is finite.
 The results of Flynn-Poonen-Schaefer \cite{flynn-poonen-schaefer} may be viewed
as some positive evidence, however, aside
from this  we do not have many other examples (however see Theorem \ref{thm:three}).
 If true, this would imply the uniform boundedness of
torsion of abelian varieties and other similar conjectures.
\end{rems}

In general the set of periodic points is of course not finite. 
However, one can often use some geometric arguments to deduce finiteness
of the number of periodic points from Theorem \ref{thm:one} as in the following:

\begin{lem}
Let $X$ be a proper variety over a finitely generated field $K$
of characteristic zero and $f:X \to X$ a morphism. Suppose that
there does not exist any positive dimensional subvariety $Y$
of $X$ such that $f$ induces an automorphism of finite order of 
$Y$. Then the number
of $f$-periodic points in $X(K)$ is finite.
\end{lem}
\begin{proof}
Theorem \ref{thm:one} implies that $f$ induces an automorphism 
of finite order on
the closure of the set of $f$-periodic points in $X(K)$.
\end{proof}

The following gives a useful method for checking the hypothesis
of the previous lemma.

\begin{lem}
Let $X$ be a projective variety over a field $K$ and $f:X \to X$
a morphism. Suppose there exists a line bundle $L$ on
$X$ such that $f^*(L) \otimes L^{-1}$ is ample. Then there is no
positive dimensional subvariety $Y$
of $X$ such that $f$ induces an automorphism of finite order of $Y$.
\end{lem}
\begin{proof}
By replacing $X$ by $f^n(X)$ for some large $n$, we may assume
that $f$ is a finite morphism. Suppose there exists a $Y$ as above
and assume that $f^m|Y$ is the identity of $Y$. Then
\[
f^m(L) \otimes L^{-1} = \bigotimes_{i = 1}^{m-1} {f^i}^*(f^*(L) \otimes L^{-1})
\]
By assumption $f^*(L) \otimes L^{-1}$ is ample so $f^m(L) \otimes L^{-1}$,
being a tensor product of ample bundles, is also ample.
But $f^m(L) \otimes L^{-1}|_Y$ is trivial, so it follows that $Y$ must be
$0$-dimensional.
\end{proof}

In  case $L$ is also ample, finiteness can also be proved using heights,
see for example \cite{kawaguchi}. One advantage of our method is that
it applies also when $f$ is an automorphism, in which case 
an ample $L$ as above can never exist. Using this, one can
extend the finiteness results of Silverman \cite{silverman-K3} to apply
to all automorphisms $f$ of projective algebraic varieties $X$ for which $1$ is
not an eigenvalue of $f^*$ on $NS(X)_{\mathbf{Q}}$.

\section{p-adic fields}

Proposition \ref{prop:one} shows that one also has boundedness of periods
for $p$-adic fields, upto powers of $p$, as long as the variety and the
morphism extend to the ring of integers. We now show that 
in fact we can bound the extra powers of $p$.

\begin{thm} Let $\mathcal{O}$ be the ring of integers in $K$, a finite
extension of $\mathbf{Q}_p$,  and let
$\mathcal{X}$ be  a proper scheme of finite type over $ Spec(\mathcal{O})$. Then
there exist a constant $M >0$  such that for any $\mathcal{O}$-morphism
$f: \mathcal{X} \to \mathcal{X}$, the periods of the $f$-periodic
 points of $\mathcal{X}(K)$ 
are all less than $M$.
\label{thm:two}
\end{thm}

If $X$ is any variety over a finite field $k$ then it is clear that
a statement analogous to the theorem holds for $X(k)$: since this is
a finite set the periods are bounded above by
$|X(k)|$, and hence are bounded independently of the morphism. 
 To bound the powers of $p$ that occur, one  sees
from the proof of Proposition 1 that it is enough to prove the following:

\begin{prop} Let $\mathcal{O}$ be the ring of integers in $K$, a 
discrete valuation ring of characteristic zero with residue
field $k$ of characteristic $p$.
 Let $(A,m)$ be a local sub-$\mathcal{O}$-algebra of
$\mathcal{O}^{p^n}$ of rank $p^n$ which
is preserved by the automorphism $\sigma$ given by cyclic permutation of the 
coordinates. Furthermore, assume that $\sigma$ acts trivially on
$m/m^2$. Then $n < r = \nu(p)$ if $p>2$
and $n \leq r$ if $p=2$, where $\nu$ is the normalized valuation
on $K$.\label{prop:two}
\end{prop}
\begin{proof}
Assume that $n \geq r$ if $p>2$ and $n > r$ if $p=2$. Since $\sigma$ acts
trivially on $m/m^2$, it follows that $\sigma^{p^t}$ acts trivially on
$m/m^{t+2}$ for all $t \geq 0$. Thus, by replacing $A$ by a quotient algebra
corresponding to the Zariski closure  in $Spec(A)$ of the  $\sigma^{p^{n-1}}$ orbit of
 any $\mathcal{O}$
valued point,
we obtain a local rank $p$ subalgebra of $\mathcal{O}^p$ which is stable
under (the new) $\sigma$ and such that $\sigma$ acts trivially on
$m/m^{r+1}$ ($m/m^{r+2}$ if $p=2$).

For $a$ in $A$ we denote by $\nu(a)$ the \emph{minimum} of the valuations of
the coordinates. Let 
\[ U(m) = \{a \in m | \nu(a_i) \neq \nu(a_j) \mbox{ for some } i,j\} \]
and let 
\[
P(m) = \{a \in U(m) | \nu(a) \leq \nu(b) \mbox{ for all } b \in U(m)\}.
\]
Suppose $\nu(a) = 1$ for $a \in P(m)$. Since $P(m) \subset U(m)$, it follows
that $\nu(\sigma^s(a) - a) = 1$ for some $s$, which in turn implies that
$\sigma^s(a) - a \notin m^2$. This is a contradiction, hence $\nu(a) > 1$
for all $a \in P(m)$.
Also, one easily sees that any element
of $m$ can be written as $a = x + b$ with $x \in \pi \cdot\mathcal{O}$ and
$b \in U(m) \cup \{0\}$, where $\pi$ in $\mathcal{O}$ is a uniformizing parameter.

Now let $a \in P(m)$ and consider $\sigma(a) -a$. Letting $t= r+1$ if $p>2$
and $t=r+2$ if $p=2$, we see that
\[
\sigma(a) - a = \sum_i \prod_{j=1}^t (x_{i,j} + b_{i,j}) ,
\]
with $x_{i,j} \in \pi \cdot\mathcal{O}$ and $b_{i,j} \in U(m) \cup \{0\}$.
Expanding the products and using the fact that  $\nu(b_{i,j}) > 1$,
 we see that
\[
\sigma(a) - a = x + \sum_k{z_kd_k} \mbox{ mod } \pi^{\nu(a)+t} ,
\]
with $x, z_k \in \pi^{t-1} \cdot\mathcal{O}$ and $d_k \in P(m)$. Further, using the
fact that $a$ is in $U(m)$, one sees that $\nu(x) = \nu(\sigma(a)-a) =\nu(a)$.
Thus, we get
\begin{equation}
0 = \sigma^p(a) - a = \sum_{i =0}^{p-1}\sigma^i(\sigma(a) -a) 
=p\cdot x + \sum_k z_k(\sum_{i=0}^{p-1} \sigma^i(d_k)) \mbox{ mod } \pi^{\nu(a)+t} .        \label{eq}
\end{equation}
Now the $d_k$'s are also in $P(m)$, so we  have 
\[
\sigma(d_k) - d_k = w_k \mbox{ mod } \pi^{\nu(a)+t-1},
\]
with $w_k$ in $\pi^t\cdot\mathcal{O}$.
This implies that
\[
\sum_{i=0}^{p-1} \sigma^i(d_k) = p \cdot d_k + \frac{(p-1)p}{2} w_k \mbox{ mod } \pi^{\nu(a)+t-1} .
\]
 Substituting this in equation
(\ref{eq}) (using that the $z_k$'s are in $\pi \cdot \mathcal{O}$) we 
get
\[
p\cdot x + \sum_k z_k (p\cdot d_k + \frac{(p-1)p}{2} w_k) = 0 \mbox{ mod } \pi^{\nu(a)+t}.
\]
We have $\nu(x) = \nu(a) = \nu(d_k) = \nu(w_k)$, $\nu(z_k) \geq t-1$ and $\nu(p)=r$.
By the choice of $t$ it follows that the only term in the above equation with valuation less than or equal to $\nu(a)+r$ is $p\cdot x$. This is a contradiction since
$t>r$.
\end{proof}

\begin{rem} The assumption of
properness is  used only to guarantee the existence of specializations. If
we consider an arbitrary separated scheme $\mathcal{X}$  of finite type
over $Spec(\mathcal{O})$, then we obtain boundedness
of the periods for the set of periodic points  in $\mathcal{X}(\mathcal{O})$. 
%One can also construct examples for which the set of periods of the
%periodic points in $\mathcal{X}(K)$ is unbounded.
\end{rem}

Theorem \ref{thm:two} can be extended to the case of some birational maps;
this can be used  to extend to the $p$-adic case
some of the results of S.~Marcello \cite{marcello} on automorphisms of affine
spaces. For $X$ a variety over
a field $K$ and $f:X \dashrightarrow X$ a rational map, we will say that
$P \in X(K)$ is \emph{f-periodic} if $f$ is defined at $P, f(P), f^2(P),\ldots $ and
$f^n(P) = P$ for some $n>0$.

\begin{thm} Let $\mathcal{O}$ be the ring of integers in $K$, a finite
extension of $\mathbf{Q}_p$,  and let
$\mathcal{X}$ be  a proper integral scheme of finite type over $ Spec(\mathcal{O})$. Then
there exist a constant $M >0$  such that for any birational map 
$f: \mathcal{X} \dashrightarrow \mathcal{X}$
with $Z(f^n) \cap Z(f^{-n}) = \emptyset$  for all $n>0$
the periods of the $f$-periodic
 points of $\mathcal{X}(K)$ 
are all less than $M$. Here for any $m \in \mathbf{Z}$, $Z(f^m)$ denotes the
exceptional locus of $f^m$.
\label{thm:three}
\end{thm}

\begin{proof}
The proof is essentially the same as that of Theorem \ref{thm:two}: Let $\mathcal{Z}$
be the Zariski closure (with the reduced induced structure) of the orbit of a 
periodic point in $P \in \mathcal{X}(K)$. We cannot conclude as before that $f$
induces an automorphism of $\mathcal{Z}$, but the assumptions imply that if $\mathcal{Z}_1$
is a connected component of $\mathcal{Z}$, then for all $n$
such that $f^n(\mathcal{Z}_1(K)) \cap \mathcal{Z}_1(K) \neq \emptyset$,
$f^n$ does induce an automorphism of $\mathcal{Z}_1$. Since the residue field of
$\mathcal{O}$ is finite, the number of connected components of $\mathcal{Z}$
is bounded independently of $f$ and $P$. This allows us to apply
Proposition \ref{prop:two} and the method of Proposition \ref{prop:one}
to conclude the proof. 
\end{proof}

\begin{rem}
In the published version of this paper it was remarked that one could use
the boundedness of periods for actions of algebraic groups to prove
boundedness of periods for all automorphisms of affine space. However this was based
on a misreading by the author of the results of Shafarevich \cite{shaf}
so that remark is not valid. 
\end{rem}

It is not difficult to construct examples of self-maps of varieties over any $p$-adic
field $K$ for which the periods are not bounded, however in all such examples known
to the author the set of $\overline{K}$-valued periodic points with a given period $n$ is
infinite for some $n$. This suggests the following:

\begin{ques}
Let $f:X\to X$ be a self-map of an algebraic variety over a $p$-adic field $K$
and suppose that the set of $\overline{K}$-valued periodic points of period $n$ is finite
for all integers $n$. Then is it true that the set of $K$-valued periodic points
is finite?
\end{ques}

A weaker question which seems likely to have a positive answer is:

\begin{ques}
Let $f:X \to X$ be a self-map  of an algebraic variety over a number field $K$
and suppose that the set of $\overline{K}$-valued periodic points of period $n$ is finite
for all integers $n$. Then is it true that the set of $K_{\wp}$-valued periodic points
is finite for all but finitely many places $\wp$ of $K$?
\end{ques}

%\vspace{4mm}

\noindent \emph{Acknowledgements.}
I would like to thank Bjorn Poonen and Ramesh Sreekantan for some interesting conversations
and helpful correspondence.

%\bibliographystyle{siam}
%\bibliography{../sources}

\end{document}